\documentclass[leqno, 11pt]{amsart}

\usepackage{color}

\usepackage[colorlinks,hyperindex, linkcolor=links,
citecolor=refs]{hyperref}

\usepackage{latexsym}
\usepackage{amsfonts}
\usepackage{amsthm}
\usepackage{amssymb}

\theoremstyle{plain}

\definecolor{links}{rgb}{0.2116,0.0104,0.7716}
\definecolor{refs}{rgb}{0.5812,0.0665,0.0659}

\newlength{\myindent}
\newlength{\mywidth}
\newcommand{\eq}{\begin{equation}}
\newcommand{\en}{\end{equation}}

\newcommand{\Rem}{\noindent{\bf Remark:}\ }

\newcommand{\Pf}{\noindent{\em Proof}:~}

\newcommand{\RR}{{\mathbb{R}}}

\newcommand{\Rn}{{\mathbb{R}}^{n}}

\newcommand{\noi}{\noindent}

\newcommand{\lap}{\Delta}

\newcommand{\dg}{\hbox{deg}}
\newcommand{\open}{\mathcal{U}}
\newcommand{\ban}{\mathcal{B}}
\newcommand{\vp}{\varphi}
\newcommand{\vpu}{\left(\varphi-u\right)}
\newcommand{\vpue}{\left(\varphi-u-\eta\right)_{_+}}
\newcommand{\haus}{\mathcal H}
\newcommand{\bdry}{\partial}
\newcommand{\bfzero}{\mathbf{0}}

\newcommand{\gt}{\rightarrow}
\newcommand{\ds}{\displaystyle}
\newcommand{\sk}{\vspace{7pt}}
\newcommand{\sm}{\vspace{4pt}}
\newcounter{theorem}
\newtheorem{lemma}[theorem]{Lemma}

\newtheorem{thm}[theorem]{Theorem}
\newtheorem{corollary}[theorem]{Corollary}
\begin{document}

\author{Alexander M. Meadows}
\date{}

\title
{
Stable and Singular Solutions of the Equation $\lap u = \frac{1}{u}$
}
\maketitle
\noindent
{\bf Abstract:} 
{\footnotesize We study properties of the semilinear elliptic equation 
$\lap u = \frac{1}{u}$ 
on domains in $\Rn$, with an eye toward nonnegative
singular  solutions as limits of 
positive smooth solutions.  We prove the nonexistence of such
solutions in low dimensions when we also require them to be stable for the
corresponding variational problem.  The problem of finding singular solutions is
related to the general study of singularities of minimal hypersurfaces of
Euclidean space.  }

\section{Introduction} \label{introsection}

One way to obtain singular minimal hypersurfaces with symmetry is 
given in the paper \cite{LMS:1} by Simon.  
Given positive $u:\Omega\subset\Rn\longrightarrow\RR$ solving the equation

\eq \label{minsurf}
{\mathcal M}u:=\sum_{i=1}^n D_{i}\left(
\frac{D_{i}u}{\sqrt{1+|Du|^2}}\right)=\frac{m}{u\sqrt{1+|Du|^2}}
\en

\noi
and an $m$-dimensional closed subgroup 
$\Gamma$ of the orthogonal group in $\RR^N$, with 
$G_{p}:=\left\{g(p):g\in\Gamma\right\}$ an orbit of maximal volume over 
$p\in S^{N-1}$, the ``symmetric graph'' $G(u)$ defined by
$$
G(u)=\left\{(x,u(x)\omega):x\in\Omega, \omega\in G_{p}\right\}\subset\RR^{n+N}
$$
\noi
will be stationary with respect to $n+m$--dimensional volume,
and will have the same regularity as the function $u$.  
Positive solutions $u>0$ are smooth.  However, 
if a sequence $u_{j}$ 
of such solutions converges continuously to a
weak solution $u\ge 0$, then $u$ will be singular exactly at the 
points where it is zero.  The corresponding
$G(u)$ will be a singular minimal submanifold. 
A degree theoretic program for obtaining such
sequences of solutions is outlined in Simon's paper.  
For more on equation~(\ref{minsurf}), see the survey 
paper by Dierkes~\cite{DU:1}.  
Our goal is to apply a similar program to
the equation
$\lap u =
\frac{1}{u}$.  The basic degree argument is presented in
Section~\ref{degreesection}.

Notice that if we linearize the left hand side of (\ref{minsurf}) in the form
$$
\sum_{i}D_{i}D_{i}u-\sum_{i,j}\frac{D_{i}uD_{j}uD_{i}D_{j}u}{1+|Du|^2}
=\frac{m}{u},
$$

\noi
the resulting equation is $\lap u=\frac{m}{u}$.  
For this equation the constant $m$ scales with the independent variable, i.e. 
$u(Cx)$ solves $\lap u = \frac{C^2 m}{u}$, so we may restrict our attention to the
case $m=1$:  
\eq \label{laplaceone}
\lap u = \frac{1}{u}
\en

Both equations~(\ref{laplaceone}) and (\ref{minsurf}), with $m=1$, 
have the particular solution
$u(x)=C|x|$ with $C=1/\sqrt{n-1}$.  By the results of
section~\ref{conicradialsection}, this solution is indeed the
limit of a sequence of positive smooth solutions. 

We note that equation~(\ref{laplaceone}) is the Euler--Lagrange equation 
for the variational integral
\eq \label{variationalintegral}
\mathcal F (u)=\int_{\Omega}\frac{1}{2}|Du|^2+\log u
\en
while equation~(\ref{minsurf}) has variational integral
$$
\int_{\Omega}u^m\sqrt{1+|Du|^2},
$$

\noi
which is the $n+m$-dimensional volume of $G(u)$.
Notice that for arbitrary positive $u$, the integral (\ref{variationalintegral})
is not bounded below for any positive boundary data $\varphi$ on $\bdry \Omega$. 
Indeed the function 
$u_\epsilon:=
\varphi +\zeta(\epsilon-\varphi)$, where $\zeta$ is a smooth cutoff function on
$\Omega$,  satisfies
$\mathcal F(u_\epsilon)\rightarrow -\infty$ as $\epsilon \rightarrow 0$.

Equation~(\ref{laplaceone}) also arises in relation to chemical catalyst kinetics
(See \cite{BN:1} and \cite{DMO:1}).  Here it is a
special case of the more general equation 
\eq \label{alphaequation}
\lap u = u^{-\alpha}, \hspace{30pt}\alpha>0
\en
For $\alpha\ne 1$, (\ref{alphaequation}) is stationary for the variational integral
\eq \label{alphaintegral}
\int_{\Omega}\frac{1}{2}|Du|^2+\frac{1}{1-\alpha}u^{1-\alpha}
\en
Note that when $0<\alpha<1$, this integral is bounded below, and by lower
semicontinuity, as in \cite{CBM:1}, minimizers among nonnegative functions
with given boundary data exist.  In fact, most results
concerning~(\ref{alphaequation}), as in
\cite{DMO:1}, are limited to the case
$0<\alpha<1$ and are results on minimizers of the variational problem.  In this
way, (\ref{laplaceone}) is an interesting limiting case.

An important difference in our method is that the singular limit must (weakly)
satisfy the PDE $\lap u = \frac{1}{u}$ on the whole domain $\Omega$.  Minimizers
of (\ref{alphaintegral}) need not satisfy the Euler--Lagrange equation
(\ref{alphaequation}) on the whole domain, in fact they need only solve the {\em
free boundary problem} on the set $\{u>0\}$.  


Despite the nonexistence of minimizers, we may consider solutions $u$ of
(\ref{laplaceone}) which are also stable for $\mathcal F(u)$ in the sense that
$\left.
\frac{d^2}{dt^2}\right|_{t=0}\mathcal F(u+t\zeta)\ge 0 $ for test functions
$\zeta$.  In this case, $u$ satisfies the stability inequality 
\eq \label{stabilityinequality}
\int_\Omega\frac{\zeta^2}{u^2}\le \int_\Omega |D\zeta|^2
\en
for all test functions $\zeta$.  We call such $u$ ``stable solutions.''  

In Sections~\ref{holdersection}, 
\ref{lowerboundsection}, and \ref{hausdorffsection}, we prove the main
results on stable solutions, namely, the H\"older continuity of stable solutions,
the nonexistence of singular stable solutions in dimension less than seven, and an
estimate on the size of the singular set of a stable solution.   The existence and
uniqueness of stable solutions is presented in Sections~\ref{generalmaximal}
and~\ref{maxstablesection}.
Here we state the main results.

\begin{thm} \label{existuniquestable}
Let $\Omega\subset\Rn$ and suppose there is a subsolution $v$ with $\lap v \ge 1/v$
on
$\Omega$ with boundary values $\varphi_0$ on $\bdry \Omega$.  Then for any
$\varphi\ge \varphi_0$, there is a unique stable solution $u$ of $\lap u =1/u$ on
$\Omega$ with boundary values $\varphi$.
\end{thm}

\begin{thm} \label{holderstableintro}
For a stable solution $u$ with boundary data $\varphi\le M$, and
for every $0<\alpha<1$ and $\tilde{\Omega}\subset\subset\Omega$, 
there is a constant $C(n,M,\alpha,\tilde{\Omega})$
such that for all $x,y\in\tilde{\Omega}$, $|u(x)-u(y)|\le C|x-y|^\alpha$.
\end{thm}

\begin{thm}
Let $2\le n\le 6$ and let $u$ be a positive stable solution of
$\lap u =
\frac{1}{u}$ on the $C^{1,1}$ domain $\Omega$, with boundary data
$\varphi\in C^{2,\alpha}(\Omega)$, $|\varphi|_{2,\alpha}\le M$,
and $\vp\ge \epsilon>0$.  Then there is a constant
$\delta=\delta(\Omega, M, \epsilon)$ such that
$u\ge \delta$ on $\Omega$. 
\end{thm}

\begin{corollary}
In dimension less than seven, there are no singular stable solutions
of~(\ref{laplaceone}).
\end{corollary}

\begin{thm}
Suppose $u$ is a limit of positive stable solutions of $\lap u = \frac{1}{u}$ on a
domain $\Omega$ with singular set $A=\{u=0\}$.  Then the Hausdorff dimension
$\hbox{dim}_{\mathcal H}(A)\le n-4-2\sqrt{2}$.
\end{thm}

\section{Basic Facts} \label{basicsection}

Positive solutions to equation~(\ref{laplaceone}) are subharmonic, 
so the maximum principle implies that they achieve their 
maximum on the boundary of $\Omega$.  
However, the difference $w=u-v$
of two solutions satisfies the equation $\lap w +\frac{w}{uv}=0$, 
and so the maximum principle 
does not guarantee that a positive solution is a unique solution of the
Dirichlet problem for its boundary data.  Indeed, solutions of the analogous
ordinary differential equation
$u^{\prime\prime}=\frac{1}{u}$ need not be unique,  and
positive radially symmetric solutions in low dimension may not be unique.

Solutions are invariant under  homothetic scaling of the graph.  
That is, if $u(x)$ is a solution on $\Omega$, then  
$u(Cx)/C$ is also a solution on the domain $\Omega/C$ for $C>0$.  In
particular, the conical solution scales under homothety of 
the graph to itself, when centered at the origin.

Nonnegative limits of positive solutions 
are singular, as in the following

\begin{lemma}
If $u\ge 0$ is a weak solution of $\lap u =\frac{1}{u}$ in a 
neighborhood of $x_0$ with $u(x_0)=0$, then $u$ is not differentiable at
$x_0$.
\end{lemma}

\Pf
Suppose we have a weak solution $u$ with $u\ge 0$, $u(x_0)=0$, and $u$
differentiable at $x_0$.  
Then $Du(x_0)=0$, and  for any $\epsilon>0$, in a small enough ball
$B=B_{\delta}(x_0)$ we have 
$0\le u\le \epsilon |x-x_0|$.  
The weak equation in this ball is
\eq \label{doubleweakeq}
\int_{B}u \lap\zeta =\int_{B}\frac{\zeta}{u}
\en

\noi
where $\zeta$ is $C^2$ with compact support in the ball.

\noi
We choose $\zeta$  radially symmetric such that $0\le \zeta\le 1$, $\zeta=1$ on 
$B_{\delta/2}$, and 
$|\lap \zeta|\le\frac{C}{\delta^2}$.
Then
$$
\int_B \frac{\zeta}{u} \ge \int_{B_{\delta/2}} \frac{1}{\epsilon|x|} = 
\frac{n\omega_{n}}{\left(n-1\right)\epsilon}\left(
\frac{\delta}{2}\right)^{n-1}
$$
where $\omega_n$ is the volume of the unit $n$-ball.  But also
$$
\int \frac{\zeta}{u} =\int_{B_{\delta}}u \lap\zeta 
 \le  
\int_{B_{\delta}}\epsilon |x| \frac{C_1}{\delta^2} = 
\frac{n\omega_n C_2 \epsilon}{n+1}\delta^{n-1}.
$$

\noi
Thus, $\epsilon^2 \ge \left(\frac{n+1}{n-1}\right)\frac{ 2^{1-n}}{C}$, 
a contradiction for 
$\epsilon$ small enough.
\sk

Using similar methods, we can derive some basic positivity results.

\begin{lemma} \label{basicpositivity}
Suppose $u$ is a positive and smooth subsolution, i.e. $\lap u \ge
\frac{1}{u}$, on a ball
$B_{2\rho}$.  Then the following hold:
\begin{enumerate}
\item $\ds \frac{1}{\rho^2} \int_{B_{2\rho}\setminus B_\rho} u^2 \ge
\omega_n \rho^n$
\item  $\ds \sup_{B_{2\rho}}u \ge
\frac{\rho}{\sqrt{2^n-1}}$
\end{enumerate}
\end{lemma}

\noi
Notice that the second property says that there do not exist solutions
on any ball that are uniformly small, and any solution
defined on all of $\Rn$ must be unbounded. 

\section{Asymptotically Conical and Radial Solutions} \label{conicradialsection}

In this section it will be convenient to consider the rescaled equation
$$
\lap u = \frac{n-1}{u}
$$
which has the conical solution $u(x)=|x|$.  We
will also use the radial variable $r=|x|$.

Caffarelli, Hardt, and Simon proved in \cite{CHS:1} the 
existence of minimal surfaces
asymptotic to minimal cones.  Following their proof we get
a similar result showing the existence of a wide variety of singular solutions of
(\ref{laplaceone})  on the ball which are asymptotic to our conical solution at the
origin.  

The wide variety of singular solutions comes from the existence of solutions with
boundary data which are small perturbations of constant boundary data equal to one
on the sphere $\bdry B_1$.  As in
\cite{CHS:1}, we are only able to specify the perturbed boundary data of the
asymptotic solution in the orthogonal complement of a finite dimensional subspace
of $L^2(\bdry B_1)$, which depends on the rate at which our solutions will
be asymptotic to the cone.  In this case, the finite dimensional subspace is the
span of the first
$J$ eigenvectors of the operator
$-\lap - (n-1)$ on the sphere, where the remaining eigenvectors have eigenvalues
$\mu$ large enough that 
$1-\frac{n}{2}+\sqrt{\frac{(n-2)^2}{4}+\mu}>m$.  The projection onto the
orthogonal complement is denoted $\Pi_J$.  The result may be expressed in terms of
the scaled H\"older norm on annuli, defined by 

$$
\left|f\right|_{2,\alpha;r}= \sum_{l=0}^2 r^l \sup_{r\le|x|\le 2r}|D^l u|+ \,
r^{k+\alpha}
\hspace{-20pt}
\sup_
{\tiny
\begin{array}{c}
x\ne y \\
r\le |x|,|y|\le 2r \\
\end{array}
} \hspace{-20pt}
\frac{|D^2f(x)-D^2f(y)|}{|x-y|^\alpha}
.
$$

\begin{thm}
Given $m>1$ and $0<\alpha<1$, there exist $\epsilon$ and $C$ depending
on $m,n$, and $\alpha$ so that for any function $g$ on $\bdry B_1$
with $|g|_{C^{2,\alpha}}<\epsilon$, there exists a solution $u$ of 
$\lap u = \frac{n-1}{u}$ on $B_1$ with $\Pi_J (u-1) =\Pi_J g$ on
$\bdry B_1$ and satisfying for $0<r<1/2$ 
$$
r^{-m}\left| u-|x|\right|_{2,\alpha;r}\le C|g|_{C^{2,\alpha}}.
$$
\end{thm}
The proof of this theorem is essentially identical to the argument given in
\cite{CHS:1}.

We call smooth positive solutions of $\lap u =
\frac{n-1}{u}$ which are radially symmetric, i.e. $u(x)=u(r)$,
$r=|x|$, ``radial solutions.''  The resulting
ordinary differential equation satisfied by
$u$ is
$u_{rr}+\frac{n-1}{r}u_r-\frac{n-1}{u}=0$, with the particular solution $u=r$.  
The following two theorems on radial solutions are the most useful for further
analysis of the PDE.  Theorem~\ref{mainradialexistence} is due to
Brauner--Nicolaenko
\cite{BN:1}, using bifurcation theory in the context of
equation~(\ref{alphaequation}).  We have alternative proofs of these and more
general facts using basic ODE techniques.

\begin{thm} \label{mainradialexistence}
For any $\epsilon>0$, solutions of the ODE problem
$u_{rr}+\frac{n-1}{r}u_r-\frac{n-1}{u}=0$ with
$u(0)=\epsilon$ and $u^\prime(0)=0$ exist uniquely on $\left[0, \infty\right)$. 
These solutions satisfy $u(r)-r= O(1)$, and consequently as $\epsilon \gt 0$, the
solutions $u(r)\gt r$ uniformly on compact subsets.
\end{thm}

\begin{thm} \label{mainradialuniqueness}
There exist constants $C_1$ and $C_2$ depending on $n$ such that on the ball
$B_1(\bfzero)$, the Dirichlet problem
$$
\begin{array}{rclcl}
\lap u &=& \frac{n-1}{u} &\hspace{20pt}&\hbox{ on }B_1 \\
u&=&C &\hspace{20pt}&\hbox{ on }\bdry B_1
\end{array}
$$
has a solution for $C>C_1$, and has a unique solution for $C>C_2$.  For $n\ge 7$,
$C_1=C_2=1$.
\end{thm}

Thus, the radial conic solution $u=|x|$ is indeed a limit of positive smooth
solutions.  An interesting further result
is that these conic solutions are stable for $n\ge 7$ and unstable
$2\le n\le 6$.  We state the result for the original equation $\lap u = 1/u$.

\begin{lemma} \label{hardyresult}
The conical solutions $u=\frac{|x|}{\sqrt{n-1}}$ are stable for $n\ge 7$ and are
unstable for $2\le n\le 6$.
\end{lemma}

\Pf
This follows from the Hardy inequality with best constant
$$
\frac{(n-2)^2}{4}\int_\Omega \frac{\zeta^2}{|x|^2}\le \int_\Omega |D\zeta|^2
$$
for all $\zeta\in C_c^1(\Omega)$.  See \cite{HLP:1}, or for a simple
proof see \cite{GAP:1}.

\section{Degree Construction} \label{degreesection}

We will use the Leray--Schauder degree with several different setups.  For the
basic theory of the degree on Banach spaces, see~\cite{KD:1}. In general we will
use the Banach space
$\ban=C^{2,\alpha}(\bar\Omega)$ and  open set $\open=\{u\in\ban: u>g,
|u|_{2,\alpha}<M_\delta\}$  where $g$ is a fixed positive bounded
function with positive minimum $\delta$.
A typical operator $T:[0,1]\times\open\rightarrow\ban$ will be defined
by $T_t(u)=v$, where $v$ is the solution of
$$
\left\{
\begin{array}{lclcl}
\lap v &=& \frac{v}{u^2} &\mbox{ on }& \Omega \\
v &=& \vp_t &\mbox{ on }& \partial \Omega
\end{array}
\right.
$$

\noi
and $\vp_t$ are boundary data continuous in $t$ with $g<\vp_t<M_\delta$.  We use
the notation
$\dg(I-T_t,\open,0)$ for the
Leray--Schauder degree invariant for fixed points of $T_t$.

In the following results, all
solutions are assumed to be positive.  Lemma~\ref{schauderbound} comes from a
basic Schauder estimate.

\begin{lemma} \label{schauderbound}
For each $0<\delta<1$, $M_\delta$ can be chosen such that any solution
$u$ of $\lap u = \frac{1}{u}$ with $\delta\le u\le \frac{1}{\delta}$
satisfies
$|u|_{2,\alpha}<M_\delta$.
\end{lemma}

\begin{lemma} \label{degreeone}
If $\vp_0\ge C_2(n)$ is constant
boundary data on the unit ball in $\Rn$, $C_2(n)$ as in
Theorem~\ref{mainradialuniqueness}, and if $\open$ is convex containing $\vp_0$
and the solution of~(\ref{laplaceone}) with data $\vp_0$, then
$\dg(I-T_0,\open,0)=1$.
\end{lemma}

\noi
{\em Proof of Lemma~\ref{degreeone}}:

\noi
By Theorem~\ref{mainradialuniqueness}, the radial solution with
$u=\vp_0$ on $\bdry B_1$ is unique.  We let $T_t(u)=v$ be the solution to
the problem
$$
\left\{
\begin{array}{lclcl}
\lap v &=& \frac{(1-t)v}{u^2} &\mbox{ on }& B_1 \\
v &=& \vp_0 &\mbox{ on }& \partial B_1
\end{array}
\right.
$$

\noi
so that the unique solution $u$ above is the unique fixed point of
$T_0$.  From $u$ we can scale to $\tilde{u}=u(\sqrt{1-t}x)$ which
satisfies $\lap \tilde{u}=\frac{1-t}{\tilde{u}}$ and is unique relative
to its boundary data, $\tilde{u}=u(\sqrt{1-t})<\vp_0$ on $\bdry B_1$.  We
may then geometrically scale $\tilde{u}$ to get $\hat{u}$ 
uniquely solving
$\lap \hat{u}=\frac{1-t}{\hat{u}}$ and $\hat{u}=\vp_0$ on $\bdry B_1$.  Note that
$\hat{u}>u$ on $B_1$.  
Thus, $T_t$ has a unique fixed point for all $t$, and in our
Leray--Schauder degree setup, there are no fixed points of $T_t$ on
the boundary of $\open$ for any convex $\open$
containing $u$ and $\vp_0$.  So, $\dg(I-T_0,\open,0)=\dg(I-T_1,\open,0)$.  But
$T_1\equiv C$.  So,
$\dg(I-T_1,\open,0)=\dg(I-C,\open,0)=\dg(I,\open,C)=1$.
\sk

\begin{lemma} \label{degreezero}
There exists $\epsilon(\Omega)$ such that if $\vp_1\le \epsilon$ 
is boundary
data on a domain $\Omega$ in $\Rn$, then no
solution with boundary data $\phi_1$ exists and
$\dg(I-T_1,\open,0)=0$.
\end{lemma}

\Pf
Lemma~\ref{degreezero} follows easily from part 2 of
Lemma~\ref{basicpositivity}





\begin{lemma} \label{finitesubsolution}
If $\Omega$ is an arbitrary domain and the function $g$ is chosen to be
the maximum of a finite collection of subsolutions, i.e. 
$\open=\{u\in C^{2,\alpha}: u>g_1, \ldots, u>g_k,
|u|_{2,\alpha}<M_\delta\}$ with 
$\lap g_k \ge
\frac{1}{g_k}$, and if $\vp_1>\max g_k$ is boundary data on $\Omega$, then
$\dg(I-T_1,\open,0)=1$.
\end{lemma}

\Pf
Consider again the map
$T_t(u)=v$ where $v$ is the solution of
$$
\left\{
\begin{array}{lclcl}
\lap v &=& \frac{tv}{u^2} &\mbox{ on }& \Omega \\
v &=& \vp_0 &\mbox{ on }& \partial \Omega
\end{array}
\right.
$$

\noi
Since $T_0$ is constant, it has degree one.  Suppose for some $t$ that
$T_t$ has a fixed point $u$ in $\bdry\open$.  Then $\lap u =\frac{t}{u}$,
and for some $g_j$, we have $u\ge g_j$ and $u(x_0)=g_j(x_0)$ for some
$x_0$.  Then, 
$$
\lap(u-g_j)\le \frac{tg_j-u}{ug_j}\le 0
\qquad
u-g_j\ge 0 \hbox{ on }\Omega.
$$

\noi
So $u-g_j$ has a zero minimum, contradicting the Hopf Maximum Principle. 
Thus, $\dg(I-T_1,\open,0)=1$.
\sk

We now outline a general method for producing ``singular sequences'' of positive
solutions to (\ref{laplaceone}) with minimum tending to zero.
In the application of the degree, let us choose $\Omega=B_1\subset\Rn$.  
Let 
$g=\delta_j>0$ so that $\open=\{u\in C^{2,\alpha}: u>\delta_j,
|u|_{2,\alpha}<M_{\delta_j}\}$, and let $\delta_j\searrow 0$. 
Using Lemmas~\ref{degreeone} and \ref{degreezero}, we may take $\vp_t$ to
be any homotopy of boundary data between $\vp_0=C\ge C(n)$ a large constant
and
$\vp_0\le \epsilon$ small.  Then since $\dg(I-T_0, \open, 0)=1$ and
$\dg(I-T_1, \open, 0)=0$, there must exist $t_j\in(0,1)$ and a fixed point
$u_j\in \bdry\open$ which solves
$$
\left\{
\begin{array}{lclcl}
\lap u_j &=& \frac{1}{u_j} &\mbox{ on }& B_1 \\
u_j &=& \vp_{t_j} &\mbox{ on }& \partial B_1
\end{array}
\right.
$$
with $\ds \min_{B_1}u_j=\delta_j$. Then the sequence $u_j$ is a ``singular
sequence" in the sense that $\min u_j\rightarrow 0$.  If $u_j\rightarrow
u$ uniformly with $u\ge 0$ and $\min u = 0$, then $u$ is a singular
solution.  Notice that if $\vp_t\in C^1$ is bounded, then we at least have
a subsequence $j^\prime$ and a $t_0$ with $\vp_{t_{j^\prime}}\rightarrow
\vp_{t_0}$.  However, we do not yet have the necessary continuity estimates
on $u_j$ to get a singular limit $u$.  

In the case $n=2$, since the $u_j$ are in particular
subharmonic, we can use the ``log trick'' (See Corollary~\ref{logtrick} below) to
show that 
$$
\int_{B_\rho} \left|Du_j\right|^2 \le \frac{C}{|\log \rho|}
$$
uniformly as $\rho\rightarrow 0$, which is just short of a modulus of continuity
estimate and also shows the $u_j$ are 
uniformly of vanishing mean oscillation as in
\cite{DPS:1}.  Even in dimension 2, subharmonicity cannot be sufficient for a
continuity estimate.  For example, the functions $u_\epsilon (x)=|x|^\epsilon$ for
$\epsilon>0$ are nonnegative, uniformly bounded on $B_1\subset \RR^2$, and
subharmonic, yet do not satisfy any continuity estimate.

\section{General Facts about Maximal Solutions} \label{generalmaximal}

A maximal solution $u$ of $\lap u = \frac{1}{u}$ satisfies 
the property that $v\le u$
for any other solution $v$ with the same boundary data as $u$.  
We show in Lemma~\ref{maximalexistence} below 
that whenever a subsolution exists for fixed boundary
data, there is also a maximal solution with that boundary data.  It
turns out that the maximal solution is also the unique stable solution.  
The existence of maximal solutions can be achieved by the usual method of
sub/supersolutions (see~\cite{DMO:1}).  We give an alternative degree method.

We will consistently use the Leray--Schauder degree with the operator
$T_t(u)=v$, where $v$ is the solution of
$$
\left\{
\begin{array}{c}
\lap v=\frac{v}{u^2} \\
v|_{\bdry\Omega}=\vp_t
\end{array}
\right. .
$$


\begin{lemma}\label{finite}
For a nonempty finite set of positive subsolutions $u_j$ with boundary
data
$\varphi_1$, there is a solution $u$ with boundary data
$\varphi_1$ such that $u\ge u_j$ for all $j$.
\end{lemma}
\sk

\Pf
Consider the open set 
$$
\open= \bigcap_j \left\{ u\in C^{2,\alpha}
: u(x)> u_j(x), |u|_{2,\alpha}<M_\delta \right\}
$$ 
in
$C^{2,\alpha}(\Omega)$, where all $u_j>\delta$ and $M_\delta$ is
chosen according to Lemma~\ref{schauderbound}.  Take data
$\varphi_0>\varphi_1$ and let
$\varphi_t$ be any smooth decreasing homotopy from 
$\varphi_0$ to $\varphi_1$.  By Lemma~\ref{finitesubsolution},
$\dg(I-T_t,\open,0)=1$ for all $t<1$, and thus there exist
solutions $u_t\in \open$ for all $t<1$.    Consider any sequence
$t_j\nearrow 1$, and corresponding solutions 
$u_{t_j}$.  Since all of these functions are uniformly bounded below, 
the Schauder estimates give us a uniform $C^3$ bound, so by 
Arzela--Ascoli, a subsequence $u_{t_k}\longrightarrow u\in \bar{\open}$
in
$C_{2,\alpha}$,  and $\lap u = 1/u$, $u|_{\bdry \Omega}=\varphi_1$, 
$u\ge u_j \:\forall j$.
\sk

\begin{lemma}\label{maximalexistence}
If there is a positive subsolution $u_0$ to $\lap u = \frac{1}{u}$ on $\Omega$ with
boundary data $\varphi_0$, and if $\varphi\ge \varphi_0$, then there is a unique
maximal solution with boundary data $\varphi$.
\end{lemma}
\sk

\Pf
Consider the collection ${\mathcal C}$ of all solutions $u\ge u_0$ with boundary
data $\varphi$.  By Lemma~\ref{finitesubsolution}, this collection is nonempty. 
$\mathcal C$ is partially ordered by the relation $u_\alpha \le u_\beta$ on
$\Omega$.  
By the Hausdorff Maximality Theorem, there exists a maximal totally
ordered subset $S$.  For any $x_0\in \Omega$ let 
$u_\alpha, u_\beta \in S$ with $u_\alpha \ne u_\beta$ and $u_\alpha\le 
u_\beta$.  By the maximum principle, we have $u_\alpha(x_0)<u_\beta(x_0)$.
Thus, $S$ can be indexed by $u_\alpha(x_0)$.  That is,
$S=\left\{u_\alpha\right\}_{\alpha\in A}$ where $\alpha=u_\alpha(x_0)$.
By the maximum principle, since $\lap u_\alpha\ge 0$, $u_\alpha\le \sup
\varphi$ for all
$\alpha$,  so $A$ is bounded above.  Let
$\alpha_\infty=\sup_A\alpha$.  We claim that
$\alpha_\infty\in A$, and $u_{\alpha_\infty}$ is a maximal solution.
Consider $u_{\alpha_j}\in S$ with $\alpha_j\nearrow\alpha_\infty$ and  
$u_{\alpha_j}\le u_{\alpha_{j+1}}$.  Since these are uniformly bounded
below and monotone increasing, the Schauder estimates and Arzela--Ascoli
give a function $u_{\alpha_\infty}$ with $u_{\alpha_j}\nearrow
u_{\alpha_\infty}$, where
$u_{\alpha_\infty}(x_0)=\alpha_\infty$ and $u_{\alpha_\infty}$ is a solution.  
For any $u_\alpha\in S$, choose $j$ large so that $\alpha_j > \alpha$.  Then 
$u_{\alpha_j}\ge u_\alpha$ by total ordering, and $u_{\alpha_\infty}\ge 
u_{\alpha_j}$ since the sequence was monotone.
Thus, $u_{\alpha_\infty}$ is an upper bound for $S$, so by maximality 
$u_{\alpha_\infty}\in S$.
To see that $u_{\alpha_\infty}$ is a maximal solution, suppose 
$v$ is another solution with $v(x_1)>
u_{\alpha_\infty}(x_1)$ for some $x_1\in\Omega$.
By Lemma~\ref{finite}, there exists a solution
$\tilde{v}$ with $\tilde{v}\ge u_{\alpha_\infty}$ and $\tilde{v}\ge v$.  
Also by the hypothesis on $x_1$, $\tilde{v}\ne u_{\alpha_\infty}$.  
But then $S\cup\{ \tilde{v}\}$ is totally ordered, contradicting
maximality.
\sk

\begin{lemma}\label{maxincrease}
If $\varphi_0<\varphi_1$, and $u_0, u_1$ are maximal solutions with
boundary data $\varphi_0$ and $\varphi_1$ respectively, then
$u_0<u_1$.
\end{lemma}
\sk

\Pf
By Lemma~\ref{finitesubsolution}, there exists a solution $u$ to $\lap u =
\frac{1}{u}$ with data $\vp_1$ and $u>u_0$, since $u_0$ is a
subsolution.  Since $u_1$ is maximal, $u_1\ge u >u_0$.
\sk

\begin{lemma} \label{maxsubdomain}
If $u$ is a maximal solution on $\Omega$ and
$\tilde\Omega\subset\Omega$ is a subdomain with continuous
boundary, then $u$ restricted to $\tilde\Omega$ is a maximal
solution with respect to its boundary data on $\bdry\tilde\Omega$.
\end{lemma}

\Pf
Suppose not.  Then there is a maximal solution $v$ on
$\tilde\Omega$ with boundary data $u$ and $v>u$ on $\tilde\Omega$. 
Let
$\mathcal U$ be the open set 
$$
\open = \left\{ w\in C^{2,\alpha}
: w>u\hbox{ on }\Omega, w>v\hbox{ on }\tilde\Omega,
|w|_{2,\alpha}<M_\delta
\right\}.
$$
Let
$\vp_1$ be any boundary data on
$\bdry\Omega$ greater than the boundary data $\vp_0$ of $u$.   As in
the proof of Lemma~\ref{finitesubsolution}, consider the operator 
$T_t(w)=\tilde w$ where $\tilde w$ is the solution of
$$
\left\{
\begin{array}{lclcl}
\lap \tilde w &=& \frac{t\tilde w}{w^2} &\mbox{ on }& \Omega \\
\tilde w &=& \vp_1 &\mbox{ on }& \partial \Omega
\end{array}
\right.
$$
Suppose $w$ is a
fixed point of $T_t$ in $\bar{\mathcal{U}}$.  By the Hopf Maximum
Principle, $w>u$ on $\Omega$.  Thus, $w>u$ on $\bdry
\tilde\Omega$.  Then we may apply the maximum principle on
$\tilde\Omega$, so $w>v$ on $\tilde\Omega$.  Thus, $w$ cannot be
on the boundary of $\mathcal U$.  So, $\dg(I-T_1,\mathcal
U,0)=\dg(I-T_0,\mathcal U,0)=1$.   Then, for any boundary data
$\vp_t>\vp_0$, there exists a solution
$w$ of $\lap w = \frac{1}{w}$ on $\Omega$ with $w=\vp_t$
on $\bdry\Omega$, $w>u$ on $\Omega$, and $w>v$ on $\tilde\Omega$. Now
we let $\vp_t$ be any smooth decreasing homotopy of boundary data
approaching $\vp_0$.  Let $w_t$ be the corresponding solutions whose
existence we just proved.  By the Schauder estimates and
Arzela--Ascoli, there exists a sequence $w_{t_j}$ with
$t_j\rightarrow 0$ such that $w_{t_j}$ converges to a solution $w$
with boundary data $\vp_0$, and with $w\ge u$ on $\Omega$ and $w\ge
v$ on $\tilde\Omega$.  Thus, $w\ge v>u$ on $\tilde\Omega$,
contradicting the maximality of $u$.
\sk


\section{Stability of Maximal Solutions} \label{maxstablesection}

Recall that stable solutions satisfy the stability
inequality~(\ref{stabilityinequality}).

\begin{lemma} \label{maxarestable}
Maximal solutions of $\lap u = 1/u$ are stable.
\end{lemma}
\sk

\Pf
Let $u_0$ be a maximal solution with data $\vp_0$, and let $\vp_t = \vp_0+t$ for
$t>0$.   By Lemma~\ref{maximalexistence}, there
exist maximal solutions $u_t$ with data $\varphi_t$ and $u_t>u_0$.  
By the Schauder estimates the $u_t$ are also bounded in $C^{4}$.
For a sequence $t_j\searrow 0$, we then have a subsequence
such that $u_{t_j}\longrightarrow \tilde{u}_0$ in $C^{2,\alpha}$, 
with $\tilde{u}_0\ge u_0$, and  
by maximality 
$\tilde{u}_0=u_0$.
%
%
Let $\delta_j=\max_\Omega (u_{t_j}-u_0)$, and let
$v_j=\frac{u_{t_j}-u_0}{\delta_j}$.  Then 
$$
\lap v_j=\frac{u_0-u_{t_j}}{u_0u_{t_j}\delta_j}.
$$
By the Schauder estimates,
$v_j$ is bounded in 
$C^{2,\alpha}$, so by Arzela--Ascoli, a subsequence
$v_j\longrightarrow v$ in
$C^2$. The function $v$ is nonnegative, not identically $0$, has
nonnegative boundary data, and satisfies the 
linearized equation $\lap v + \frac{1}{u_0^2}v=0$.  By the maximum principle,
$v>0$ in $\Omega$. The weak equation for $v$ is then
$$
\int_\Omega Dv\cdot D\zeta = \int_\Omega \frac{\zeta v}{u_0^2}.
$$
We use the test function $\frac{\zeta^2}{v}$ for $\zeta$, and  
Cauchy--Schwartz to get the desired inequality 
$$
\int\frac{\zeta^2}{u_0^2}\le\int|D\zeta|^2
$$
\noindent
for all compactly supported $\zeta$.
\sk

In fact, the maximal solution for given boundary data is the only stable solution. 

\begin{lemma} \label{maxunique}
The maximal solution for given boundary data $\varphi$ is the unique stable
solution with data $\varphi$.
\end{lemma}

\Pf
Let $u$ be the maximal solution and let $v$ be any other positive solution.  Then
$u=v+w$ where $w>0$ in $\Omega$ and $w=0$ on $\bdry\Omega$.  Thus, 
$$
\lap w = \frac{1}{v+w}-\frac{1}{v} = \frac{-w}{v(v+w)}>\frac{-w}{v^2}
$$
and so, integrating by parts with $w$, 

$$
\int_\Omega \frac{w^2}{v^2} >  \int_\Omega w\lap w =\int_\Omega |Dw|^2
$$

\noi
and $v$ cannot be stable.

Theorem~\ref{existuniquestable} now follows from Lemmas~\ref{maximalexistence},
\ref{maxarestable}, and~\ref{maxunique}

\section{H\"older Continuity of Stable Solutions} \label{holdersection}

We now prove a main result that stable solutions are locally uniformly H\"older
continuous.

\begin{thm} \label{holderstable}
For a stable solution $u$ with boundary data $\varphi\le M$, and
for every $0<\alpha<1$ and $\tilde{\Omega}\subset\subset\Omega$, 
there is a constant $C(n,M,\alpha,\tilde{\Omega})$
such that for all $x,y\in\tilde{\Omega}$, $|u(x)-u(y)|\le C|x-y|^\alpha$.
\end{thm}
\sm

\Pf
Let $u$ be a smooth positive stable solution on $\Omega$. 
The weak form of the equation is 
\eq \label{weak}
\int Du\cdot D\zeta =-\int\frac{\zeta}{u}
\end{equation}

\noi
for $\zeta$ compactly supported.  Here and throughout the rest of the
proof, all integrals are taken over the domain $\Omega$. Substituting
$u\zeta^2$ for
$\zeta$ in (\ref{weak}) yields

$$
\int |Du|^2\zeta^2 + 2\int u\zeta Du\cdot D\zeta = -\int \zeta^2
$$
and applying the Cauchy--Schwartz inequality, we get
\eq\label{dirichletintegral}
\int \zeta^2\left(\frac{1}{2}|Du|^2+1\right)\le 2\int u^2|D\zeta|^2 .
\end{equation}

\noi
Substituting $\frac{\zeta^2}{u}$ in (\ref{weak}) and again using
Cauchy--Schwartz gives

\eq\label{oneonu}
\int\frac{|Du|^2}{u^2}\zeta^2\le 2\int\frac{\zeta^2}{u^2} 
+4\int |D\zeta|^2 .
\end{equation}

\noi
Differentiating the equation with respect to $x_l$, and using subscripts to
denote differentiation, we have
$\lap u_l=-\frac{u_l}{u^2}$, and thus the weak equation
$$
\int u_{li}\zeta_i=\int\frac{u_l}{u^2}\zeta,
$$  
where we sum on the repeated index $i$.  Note that this equation is
equivalent to using
$\zeta_l$ in (\ref{weak}) and integrating by parts.  
Substituting
$u_l\zeta^2$ for $\zeta$ gives
$$
\int u_{li}u_{li}\zeta^2 + 2\int u_{li}u_l \zeta \zeta_i =
\int\frac{u_l^2}{u^2}\zeta^2,
$$
or, after summing on $l$,
\eq \label{differentiated}
\int |D^2u|^2\zeta^2+2\int
u_lu_{li}\zeta\zeta_i=
\int\frac{|Du|^2}{u^2}\zeta^2
\en

\noi
It is now convenient to use the variable $v=\sqrt{1+|Du|^2}$, where

\noi
$$
v_i=\frac{1}{v} u_j u_{ji}
, \qquad
|v_i|^2\le\sum_j|u_{ij}|^2,
$$ 
and
$$
|Dv|^2\le\sum_{ij}|u_{ij}|^2=|D^2u|^2 .
$$

\noi
Also, $vv_i=u_lu_{li}$.
Replacing in equation~(\ref{differentiated}) gives
\eq\label{newdiff}
\int |Dv|^2\zeta^2+2\int vv_i\zeta\zeta_i \le
\int\frac{|Du|^2}{u^2}\zeta^2 .
\end{equation}

\noi
If $u$ is stable, it additionally satisfies (\ref{stabilityinequality}), 
$$
\int\frac{\zeta^2}{u^2}\le\int |D\zeta|^2 .
$$

\noi
Now we can combine our inequalities:
$$
\begin{array}{rll}
\int |Dv|^2\zeta^2+2\int vv_i\zeta\zeta_i &\le
\int\frac{|Du|^2}{u^2}\zeta^2 \hspace{10pt} &\hbox{ by (\ref{newdiff})}
\\ &\le 2\int\frac{\zeta^2}{u^2} 
+4\int |D\zeta|^2 \hspace{10pt} &\hbox{ by (\ref{oneonu})} \\
&\le 6\int|D\zeta|^2  \hspace{10pt} &\hbox{ by (\ref{stabilityinequality})}
\end{array}
$$

\noi
and we get the main inequality for stable
solutions:
\eq\label{mainstable}
\int |Dv|^2\zeta^2 +2\int vv_i\zeta\zeta_i \le
6\int |D\zeta|^2
\end{equation}

\noi
We use this and the Sobolev Inequality to iteratively
estimate integrals $\int v^q \zeta^\beta$ for $q\ge 0$.
In fact, the estimate for $q=2$ and $\beta=2$ is contained
in~(\ref{dirichletintegral}).
For $q>2$, replace $\zeta$ in~(\ref{mainstable})
by $v^q\zeta$ to get
\begin{eqnarray*}
\int|Dv|^2v^{2q}\zeta^2&+&2q\int |Dv|^2v^{2q}\zeta^2 
+ 2\int v^{2q+1}\zeta v_i\zeta_i \le \\
&\le& 
6\int \left| qv^{q-1}\zeta Dv + 
v^{q}D\zeta\right|^2
\end{eqnarray*}
and so
\begin{eqnarray}
\int|Dv|^2v^{2q}\zeta^2&+&2q\int |Dv|^2v^{2q}\zeta^2 \le 
\label{dagger}
\\ 2\int |Dv||D\zeta|v^{2q+1}\zeta&+& 
12q^2\int |Dv|^2v^{2q-2}\zeta^2 + 
12\int v^{2q}|D\zeta|^2 \nonumber
\end{eqnarray}

\noi
where we have used the squared triangle inequality $(a+b)^2 \le
2a^2+2b^2$.   We use the Cauchy--Schwartz inequality to eliminate the
second term on the top line of~(\ref{dagger}) with the first term on
the bottom line.
\begin{eqnarray}
\int |Dv|^2v^{2q}\zeta^2 &\le& 
\frac{1}{2q}\int v^{2q+2}|D\zeta|^2 + 
12\int v^{2q}|D\zeta|^2 + \label{bootstrapiterate} \\
& & + \hspace{2pt} 12q^2\int |Dv|^2v^{2q-2}\zeta^2 \nonumber
\end{eqnarray}

\noi
Notice at this point that the last term on the right side of the
inequality is the same as the left hand side with a lower power of
$v$.  So, we can apply~(\ref{bootstrapiterate}) to that term
iteratively until
$2q-2$ is less than zero, and use the fact that $v\ge 1$ to get

\eq
\int |Dv|^2v^{2q}\zeta^2 \le
C(n,q)\int v^{2q+2}|D\zeta|^2 .
\end{equation}

\noi
We rewrite the equation above as
\eq
\int |D(v^{q+1}\zeta)|^2\le C(n,q)\int v^{2q+2}|D\zeta|^2 ,
\en

\noi
replace $q$ by $q-1$, and apply the Sobolev inequality 
to get
\eq
\left(\int v^{2q\kappa}\zeta^{2\kappa}\right)^{\frac{1}{\kappa}} \le
C(n,q)\int v^{2q}|D\zeta|^2\hspace{10pt} \hbox{ for }q\ge 1
\end{equation}

\noi
with $\kappa=\frac{n}{n-2}$ or 
$\kappa=2$ if $n=2$.
We now replace $\zeta$ by $\zeta^\beta$, 
and we fix $\zeta$ so that
$|D\zeta|^2$ is bounded pointwise by $C(\tilde\Omega)$ and $\zeta=1$ on
$\tilde{\Omega}$.  

\eq\label{iteratable}
\left(\int v^{2q\kappa}\zeta^{2\kappa\beta}\right)^{\frac{1}{2q\kappa}} 
\le C\left(\int v^{2q}\zeta^{2\beta-2}\right)^{\frac{1}{2q}}
\end{equation}

\noi
where the constant $C$ now depends
on $n$, $q$, $\tilde{\Omega}$, and $\beta$.
Now we iterate the inequality~(\ref{iteratable}) with 
$q=1,\kappa, \kappa^2, \ldots$ and corresponding
$\beta=2$, $\beta_1, \beta_2, \ldots$, where  
$\beta_j=1+\kappa+\kappa^2+\cdots+2\kappa^j$.  Then we have
\eq
\left(\int
v^{2\kappa^m}\zeta^{2(2\kappa^m+\kappa^{m-1}+\cdots+\kappa)}
\right)^{\frac{1}{2\kappa^m}} 
\le C(n,m,\tilde{\Omega})\left(\int v^{2}\zeta^2 \right)^{\frac{1}{2}} .
\end{equation}

\noi
By~(\ref{dirichletintegral}), 

\begin{eqnarray*}
\left(\int
v^{2\kappa^m}\zeta^{2(2\kappa^m+\kappa^{m-1}+\cdots+\kappa)}
\right)^{\frac{1}{2\kappa^m}} 
&\le & C(n,m,\tilde{\Omega})\left(\int u^2 \right)^{\frac{1}{2}} \\
&\le& 
C(n,m,\tilde\Omega, M)
\end{eqnarray*}

\noi
So, on the subdomain $\tilde\Omega$ we now have a bound for the Sobolev norm:
\eq\label{finalstable}
\|u\|_{W^{1,2\kappa^m}(\tilde\Omega)}\le
C(n,m,\tilde\Omega,M)
\end{equation}

\noi
By the Sobolev Imbedding Theorem, to each $\alpha$ in the statement of the
theorem there corresponds an $m$ in equation~(\ref{finalstable}) depending
on $n$ and $\alpha$ such that we have a bound on
$|u|_{C^{0,\alpha}(\tilde\Omega)}$. This completes the proof.
\sk

\Rem This theorem may be just short of a sharp interior regularity estimate, since
the known conical example solutions are at worst Lipschitz.  For
equation~(\ref{alphaequation}) with $0<\alpha<1$, a sharp estimate for
solutions of the free boundary problem minimizing the variational integral
was given by Phillips in~\cite{P:1}.

\section{Lower Bounds for Stable Solutions in Low
Dimensions}\label{lowerboundsection}

Recall that for $n\ge 7$, the radial solutions are
unique for their Dirichlet boundary data, therefore maximal and
stable.  In particular, the conical solution is a stable singular
solution for $n\ge 7$.  
However, for $2\le n\le 6$,
the conical solution does not satisfy the stability
inequality~(\ref{stabilityinequality}) by Lemma~\ref{hardyresult} and is not
maximal for its boundary data.  Thus it cannot be the limit of stable radial
solutions and, for $2\le n\le 6$, the stable radial solutions are bounded below
by a constant.  The next result generalizes this lower bound to all stable
solutions in a compact subdomain.  Theorem~\ref{totallowerbound} gives the same
result on the entire domain.

\begin{thm} \label{lowerboundall}
Let $2\le n\le 6$ and let $u$ be a positive stable solution of
$\lap u =
\frac{1}{u}$ on the domain $\Omega$ with $u\le M$, and let
$\tilde\Omega\subset\subset\Omega$ be a compact subdomain. 
Then there is a constant $\delta=\delta(n,\tilde\Omega,M)>0$ such that $u\ge
\delta$ on $\tilde\Omega$.  
\end{thm}

\Pf
The proof follows from the estimate 
\eq \label{pintegral}
\int_{\tilde\Omega} u^{-p}\le C(\tilde\Omega, p)\hspace{10pt} \hbox{
for }p<4+2\sqrt{2}.
\en

\noi
Notice that the restriction on $p$ allows for $p\ge n$ as long as $n\le
6$.  For $u>0$ smooth, we use the stability
inequality~(\ref{stabilityinequality}) with  the test function 
$\zeta u^{-q}$.  Then for $\epsilon>0$, 
\begin{eqnarray}
\int u^{-2q-2}\zeta^2 \hspace{-5pt} &\le & \hspace{-5pt}\int |u^{-q}D\zeta -q
u^{-q-1}\zeta Du|^2 \nonumber \\
& \le & \hspace{-5pt}\int u^{-2q}|D\zeta|^2 + 2|q| u^{-2q-1}\zeta|Du||D\zeta|+q^2
u^{-2q-2}\zeta^2|Du|^2 \nonumber \\
& \le & \hspace{-5pt}\left( 1+ \frac{|q|}{2\epsilon}\right) \int u^{-2q}|D\zeta|^2
+ 
\left( q^2+ 2|q|\epsilon \right)\int u^{-2q-2}\zeta^2 |Du|^2
 \label{lowerstab}
\end{eqnarray}

\noi
Again, all integrals in this proof are taken over the domain
$\Omega$.  Notice that in every integral the integrand has compact
support in $\Omega$.  We will also use the weak form of the
equation~(\ref{weak})  with the test function
$\zeta^2 u^{-\beta}$, $\beta>0$ to get
\eq \label{weakdu}
\beta\int u^{-\beta-1}\zeta^2|Du|^2  \le  
\int u^{-\beta-1}\zeta^2+2\int
u^{-\beta}\zeta|Du||D\zeta|
\en
and using Cauchy-Schwartz, for any $\delta>0$, 
\eq \label{weakdu2}
\left(\beta-2\delta\right)\int u^{-\beta-1}\zeta^2|Du|^2 \le 
\int u^{-\beta-1}\zeta^2 + \frac{1}{2\delta}\int u^{-\beta+1}|D\zeta|^2
\en

\noi
Replacing $\beta$ by $2q+1$ and combining with~(\ref{lowerstab}), we get
%
for $q>\frac{-1}{2}$ and $\epsilon, \delta>0$, 

\begin{eqnarray*}
\int u^{-2q-2}\zeta^2 & \le & \left( 1+\frac{|q|}{2\epsilon} 
+ \frac{q^2 + 2|q|\epsilon}{2\delta\left(2q+1-2\delta\right)} \right)
\int u^{-2q}|D\zeta|^2 + \\
& + & 
\left( \frac{q^2+ 2|q|\epsilon}{2q+1-2\delta}\right) 
\int u^{-2q-2}\zeta^2
\end{eqnarray*}

\noi
Then for $1-\sqrt{2}<q<1+\sqrt{2}$ and $\epsilon$ and $\delta$ small 
enough depending on $q$, the coefficient in the last term above is less 
than one. So, 

\eq 
\int\zeta^2 u^{-2q-2} \le C(q)\int u^{-2q}|D\zeta|^2
\en

\noi
and now assuming $q>0$, we replace $\zeta$ by $\zeta^{q+1}$:
\eq \label{utoqinequality}
\int \left(\frac{\zeta}{u}\right)^{2q+2} \le C(q)\int 
\left(\frac{\zeta}{u}\right)^{2q} |D\zeta|^2
\en

\noi
Now we use Young's inequality in the form 
$$
ab \le \frac{\epsilon^\alpha}{\alpha}a^\alpha + 
\frac{\alpha-1}{\alpha \epsilon^{\frac{\alpha}{\alpha-1}}}
b^{\frac{\alpha}{\alpha-1}}
$$
with $\alpha=\frac{q+1}{q}$ and replace $q$ by $\frac{p}{2}-1$.  Then 
for $2<p<4+2\sqrt{2}$,
\eq \label{pintegralzetabound}
\int \left( \frac{\zeta}{u} \right)^{p} \le C(p)\int 
 |D\zeta|^p
\en

\noi
So we get
equation~(\ref{pintegral}) for any $p<4+2\sqrt{2}$.  We now recall
our continuity estimate for stable solutions, that for any
$\alpha<1$, the H\"older norm $|u|_{0,\alpha,\tilde\Omega}\le 
C(\tilde\Omega,\alpha, M)$.

\noi
Now let $\tilde\Omega\subset\hat\Omega\subset\subset\Omega$ with
$\hbox{dist}(\tilde\Omega,\bdry\hat\Omega)>\rho$, so that for any
$x\in\tilde\Omega$, $B_\rho(x)\subset\hat\Omega$.  For
$x_0\in\tilde\Omega$, let $r=|x-x_0|$ and suppose
$u(x_0)=\epsilon$.  Then for $x\in B_\rho(x_0)\subset\hat\Omega$, 
$u(x)\le \epsilon+C(\hat\Omega, \alpha, M)r^\alpha$.  Then
$$
\int_{\hat\Omega}u^{-p}\ge \hspace{-1pt} \int_{B_\rho(x_0)}\hspace{-7pt}u^{-p} \ge
\hspace{-1pt}\int_{B_\rho(x_0)}\left(\epsilon+Cr^{\alpha}\right)^{-p} \ge
n\omega_n\int_0^\rho \left(\epsilon+Cr^{\alpha}\right)^{-p}r^{n-1}
dr
$$
 
\noi
We may choose $\alpha$ and $p$ large enough that $n-1-\alpha p<-1$
for $n\le 6$.  Then for $\epsilon$ small enough, we have a
contradiction of equation~(\ref{pintegral}) with $\hat\Omega$ in
place of $\tilde\Omega$.  This completes the theorem.

\sk

We note that we did not need to use the continuity estimate in the
above proof.  In fact, equation~(\ref{pintegral}) together with the
$L^p$ estimates gives an estimate for $u\in
C^{1,1-\frac{n}{p}}(\tilde\Omega)$ for $n\le 6$ and
$p<4+2\sqrt{2}$.  We use a similar method below to get a lower bound on the whole
domain.  First we present an interesting corollary.

\begin{corollary} \label{logtrick}
There are no complete stable solutions of $\lap u = \frac{1}{u}$ on all of $\Rn$
for $2\le n\le 6$.
\end{corollary}

\Pf
Suppose not.  
From (\ref{pintegralzetabound}), for $2\le n\le 6$ we have
$$
\int \left(\frac{\zeta}{u}\right)^n\le C(n)\int \left|D\zeta\right|^n
$$
We choose $\zeta$ equal to one on the ball $B_R$, equal to zero outside $B_{R^2}$,
and equal to $2-\frac{\log |x|}{\log R}$ on $B_{R^2}\setminus B_R$.  Then, using
the variable $r=|x|$, we have
$$
\int_{B_R}\frac{1}{u^n}\le C\int_R^{R^2} \frac{r^{n-1}}{r^n(\log R)^n} dr 
\le \frac{C}{(\log R)^{n-1}}
$$
and the result follows letting $R\rightarrow \infty$.  We thank Neshan
Wickramasekera for pointing out this trick, which also appeared in reference to
the Bernstein Theorem in \cite{LMS:4}.

\begin{thm} \label{totallowerbound}
Let $2\le n\le 6$ and let $u$ be a positive stable solution of
$\lap u =
\frac{1}{u}$ on the $C^{1,1}$ domain $\Omega$, with boundary data
$\varphi\in C^{2,\alpha}(\Omega)$, $|\varphi|_{2,\alpha}\le M$,
and $\vp\ge \epsilon>0$.  Then there is a constant
$\delta=\delta(\Omega, M, \epsilon)$ such that
$u\ge \delta$ on $\Omega$.  
\end{thm}

The theorem follows from the following lemma.

\begin{lemma} \label{totallowerboundestimate}
Let $2\le n\le 6$ and let $u$ be a positive stable solution of
$\lap u =
\frac{1}{u}$ on the domain $\Omega$, with Lipschitz boundary data
$\varphi$, 
and $\vp\ge \epsilon>0$.  Let $2\le p < 4 + 2\sqrt{2}$.  Then there
is a constant $C(p, |D\vp|)$ such that 
$$
\int_\Omega \frac{1}{u^p} \le
\frac{C|\Omega|}{\epsilon^p}.
$$
\end{lemma}

\Rem
Notice that in Lemma~\ref{totallowerboundestimate} 
there is no assumption on the smoothness of
the domain.
\sm

\Pf
Let $\epsilon>0$ as in the statement and assume $u>0$ is a stable solution.  
Let $\eta >0$ and first consider the stability 
inequality with the test function $\zeta=\vpue$:
\eq \label{vpuestability} 
\int \frac{\vpue^2}{u^2} \le \int |D \vpue |^2
\en
So, 
\begin{eqnarray*}
\int |D \vpue |^2 \hspace{-8pt}&=& \hspace{-8pt} \int\hspace{-3pt} D\varphi
\cdot
\hspace{-2pt}D\vpue \hspace{-1pt}- \int \hspace{-3pt}Du\cdot\hspace{-2pt} D\vpue
\\ 
\hspace{-8pt}&=&\hspace{-8pt} \int D\varphi \cdot \hspace{-2pt}D\vpue +
\int\frac{\vpue}{u} 
\end{eqnarray*}
\begin{eqnarray*}
 & \le &  \int |D\varphi|^2 + \frac{1}{4}\int |D\vpue|^2 + 
\frac{1}{2}\int \frac{\vpue^2}{u^2} + \frac{1}{2}\int 1 \\
& \le &  \int |D\varphi|^2 + \frac{3}{4}\int |D\vpue|^2 + \frac{1}{2}\int 1 \hspace{20pt}\hbox{ (by~(\ref{vpuestability}))}
\end{eqnarray*}
So, 
\eq \label{dvpueinequality}
\int |D \vpue |^2 \le C\int \left(1+|D\varphi|^2\right)
\en

\noi
Now we use~(\ref{utoqinequality}) with the same test function, and replace 
$q$ by $p/2-1$ so that for $2<p<4+2\sqrt{2}$, 

\eq \label{mainpbound}
\int \left(\frac{\vpue}{u}\right)^p \le 
C(p) \int \left(\frac{\vpue}{u}\right)^{p-2}|D\vpu|^2.
\en

\noi
Recall from equation~(\ref{weakdu2}) that for $\beta>1$,
$$
\int u^{-\beta}\zeta^2|Du|^2 \le C\int u^{-\beta}\zeta^2 + C \int u^{2-\beta}|D\zeta|^2.
$$
So, assuming $p>4$, we replace $\beta$ by $p-2$ and $\zeta$ by $\vpue^{\frac{p-2}{2}}$
and we have 
\begin{eqnarray*}
\int \hspace{-2pt}\left(\frac{\vpue}{u}\right)^{\hspace{-1pt}p-2}
\hspace{-12pt}|Du|^2
\hspace{-5pt}&\le&
\hspace{-5pt}C\int\hspace{-2pt}\left(\frac{\vpue}{u}\right)^{p-2}  \\ 
&+& \hspace{-5pt}C\int \hspace{-2pt}
\left(\frac{\vpue}{u}\right)^{p-4}\hspace{-10pt}|D\vpue|^2
\end{eqnarray*}

\noi
We use this with~(\ref{mainpbound}) to get

\begin{eqnarray*}
\int \left(\frac{\vpue}{u}\right)^p + 
\int \left(\frac{\vpue}{u}\right)^{p-2}|D\vpu|^2 \le \\
C\int \left(\frac{\vpue}{u}\right)^{p-2} + 
C\int \left(\frac{\vpue}{u}\right)^{p-4}|D\vpu|^2
\end{eqnarray*}
where the constant $C$ now depends also on $|D\vp|$.
Now we can apply Young's inequality
twice. Then

$$
\int \frac{\vpue^p}{u^p}
+
\frac{\vpue^{p-2}}{u^{p-2}}
|D\vpu|^2
\hspace{-1pt}\le 
\hspace{-1pt}C\hspace{-1pt}\int 1 + 
|D\vpu|^2
$$

\noi
We then use our estimate~(\ref{dvpueinequality}), and let
$\eta$ tend to zero.

$$
\int \left(\frac{\vpu_{_+}}{u}\right)^p \le C|\Omega|
$$

\noi
Now for $\vp>2\epsilon$, 

$$
\frac{1}{u^p}\le
\frac{1}{\epsilon^p}\left(\left(\frac{\vpu_{_+}}{u}\right)^p+1\right)
$$
so
$$
\int \frac{1}{u^p} \le \frac{(C+1)|\Omega|}{\epsilon^p}
$$
as required.
\sk

\Rem
Notice that for $p=2$,we get a stronger result y removing the dependence ov $C$ on $|D\varphi|$ and using inequalities~(\ref{vpuestability}) and (\ref{dvpueinequality}).  Namely, 
\eq
\int_\Omega \frac{1}{u^2}\le \frac{C}{\epsilon^2}\int_\Omega
\left(1+|D\vp|^2\right)
\en
where we may assume that $\vp$ is merely in $W^{1,2}(\Omega)$.
\sk 

\noi
{\em Proof of Theorem~\ref{totallowerbound}:}
The lemma demonstrates the inequality 

\eq \label{pintegralbound}
\int \frac{1}{u^p}\le C(p, \Omega, \epsilon, |D\vp|)
\en

\noi
But of course $\frac{1}{u^p}=(\lap u)^p$ by the equation, and we
can apply the $L^p$ estimates (see \cite{GT:1} 9.14) related to the Calderon--Zygmund
Inequality.  So,
$$ 
\left\| u \right\|_{W^{2,p}} \le C(p,\Omega, M, \epsilon)
$$
and by the extended Sobolev Embedding Theorem, 

$$
\left| u \right|_{C^{1,1-n/p}} \le C(p,\Omega, M, \epsilon)
$$
which in particular implies a uniform Lipschitz bound on $u$.  Then 
if $u$ achieves the value $\delta$ at a point $x_0\in\Omega$,
$$
\int \frac{1}{u^p} \ge \int \frac{1}{(\delta + Cr)^p}
$$
where $r=x-x_0$, a contradiction of (\ref{pintegralbound}) for
$\delta<\delta(p,\Omega, M, \epsilon)$ and $p>n$.
\sk

\Rem
With this lower bound we in fact have complete regularity of stable
solutions for $n\le 6$.  From $u\in C^{1, \alpha}$ and thus (by the
lower bound) $\frac{1}{u}\in C^{1, \alpha}$, we can apply Holder
estimates to get continuous derivatives of all orders on the interior
of the domain.

\section{Hausdorff Dimension of Singular Sets of Stable
Solutions}\label{hausdorffsection}

We use Hausdorff dimension as described in~\cite{LMS:2}.

\begin{thm}
Suppose $u$ is a limit of positive stable solutions of $\lap u = \frac{1}{u}$ on a
domain $\Omega$ with singular set $A=\{u=0\}$.  Then the Hausdorff dimension
$\hbox{dim}_{\mathcal H}(A)\le n-4-2\sqrt{2}$.
\end{thm}

\Pf
We will show for any ball $B_\rho$ of radius $\rho$ whose closure is contained in $\Omega$, and any
$\beta>n-4-2\sqrt{2}$, that the Hausdorff Measure $\haus^{\beta}\left(A\cap
B_{\rho/2}\right)<\infty$.  First, for any $\delta$ with $0<\delta<\rho/4$, we cover $A\cap B_{\rho/2}$ 
by cubes
$Q_j$ of side length $2\delta$ with disjoint interiors, $j=1, \ldots, N$.  Let
$p<4+2\sqrt{2}$.  By (\ref{pintegral}), we have 
$$
\int_{\bigcup Q_j}\frac{1}{u^p}\le \int_{B_{\rho/2}}\frac{1}{u^p}\le
K<\infty
$$
with $K$ independent of $\delta$.  
By Theorem~\ref{holderstable}, for any $0<\alpha<1$, we have $u(x)\le C
\left(\hbox{dist}(x,A)\right)^{\alpha}$.  Thus, assuming all the $Q_j$ intersect $A$,  
$$
\int_{Q_j}\frac{C}{u^p} \ge
\int_{Q_j}\frac{1}{\left(\hbox{dist}(x,A)\right)^{\alpha p}} \ge
\int_{0<x_i<2\delta} \frac{1}{\left(x_1^2+ \cdots x_n^2\right)^{\alpha p/2}} \ge 
$$
$$
\ge
\frac{1}{2^n}\int_0^\delta \frac{n\omega_n r^{n-1}}{r^{\alpha p}}\, dr \ge 
\frac{n\omega_n}{2^n (n-\alpha p)}\delta^{n-\alpha p} 
$$
Choose $\alpha$ so that $\beta = n-\alpha p$.  Then we have
$$
\haus^\beta_\delta\left(A\cap B_{\rho/2}\right) \le C^\prime\sum_{Q_j}\delta^\beta  \le C^{\prime\prime}K<\infty
$$

independent of $\delta$, and the result follows.
\sk

\Rem
Recall that the analogous equation~(\ref{alphaequation}) with $0<\alpha<1$ has 
actual {\em minimizers} for the variational problem.  In~\cite{P:2}, Phillips proved a 
Hausdorff estimate on the {\em free boundary} for minimizers.  Note that in the free
boundary problem, the solution is allowed to vanish completely and not satisfy the PDE 
on an interior set of positive measure.  The technique above gives an estimate on the 
size of the singular set of a {\em singular solution} which is a limit of positive solutions 
satisfying the PDE on the whole domain.  

\sm

\Rem
The major results of this paper extend to the equation~(\ref{alphaequation})  with
$0<\alpha<1$, except for two points.  Solutions of (\ref{alphaequation}) which
achieve the value zero need not be singular, and the results of sections 8 and 9
are more complicated, with dimensions depending on $\alpha$.
\sm

\noindent
{\bf Acknowledgement: }
Part of this work is contained in the author's doctoral dissertation.  He would
like to thank his advisor, Professor Leon Simon.

    \bibliographystyle{plain}
    \bibliography{mybib}

\vspace{12pt}

\noi
Department of Mathematics

\noi
Cornell University

\noi
Ithaca, NY 14853

\noi
email: {\tt meadows@math.cornell.edu}

\vfill
\end{document}